# A CLT FOR REGULARIZED SAMPLE COVARIANCE MATRICES

By Greg W. Anderson and Ofer Zeitouni[1]

*University of Minnesota*

We consider the spectral properties of a class of *regularized estimators* of (large) empirical covariance matrices corresponding to stationary (but not necessarily Gaussian) sequences, obtained by *banding*. We prove a law of large numbers (similar to that proved in the Gaussian case by Bickel and Levina), which implies that the spectrum of a banded empirical covariance matrix is an efficient estimator. Our main result is a central limit theorem in the same regime, which to our knowledge is new, even in the Gaussian setup.

**1. Introduction.** We consider in this paper the spectral properties of a class of *regularized estimators* of (large) covariance matrices. More precisely, let $X = X^{(p)}$ be a data matrix of $n$ independent rows, with each row being a sample of length $p$ from a mean zero stationary sequence $\{Z_j\}$ whose covariance sequence satisfies appropriate regularity conditions (for details on those, see Assumption 2.2). Let $X^T X$ denote the empirical covariance matrix associated with the data. We recall that such empirical matrices, as well as their centered versions $(X - \tilde{X})^T (X - \tilde{X})$, where $\tilde{X}_{ij} = n^{-1} \sum_{k=1}^n X_{kj}$, are often used as estimators of the covariance matrix of the sequence $\{Z_j\}$, see [2]. We remark that the information contained in the eigenvalues of the covariance matrix is often of interest, for example, in principal component analysis or applications in signal processing.

In the situation where both $p$ and $n$ tend to infinity, it is a standard consequence of random matrix theory that these estimators may not be consistent. To address this issue, modifications have been proposed (see [3]) to which we refer for motivation, background and further references. Following the approach of [3], we consider regularization by *banding*, that is, by replacing those entries of $X^T X$ that are at distance exceeding $b = b(p)$ away from the diagonal by 0. Let $Y = Y^{(p)}$ denote the thus regularized empirical matrix.

Received January 2007; revised May 2007.
[1]Supported in part by NSF Grant DMS-05-03775.
*AMS 2000 subject classifications.* Primary 62H12; secondary 15A52.
*Key words and phrases.* Random matrices, sample covariance, regularization.

We focus on the empirical measure of the eigenvalues of the matrix $Y$. In the situation where $n \to \infty$, $p \to \infty$, $b \to \infty$ and $b/n \to 0$ with $b \leq p$, we give in Theorem 2.3 a law of large numbers (showing that the empirical measure can be used to construct an efficient estimator of averages across frequency of powers of the spectral density of the stationary sequence $\{Z_j\}$), and in Theorem 2.4, we provide a central limit theorem for traces of powers of $Y$. We defer to Section 9 comments on possible extensions of our approach, as well as on its limitations. We note that in the particular case of Gaussian data matrices with explicit decay rate of the covariance sequence, and further assuming $b \sim (\sqrt{n}/\log p)^\alpha$ for some constant $\alpha > 0$, the law of large numbers is contained (among many other things) in [3], Theorem 1. But even in that case, to our knowledge, our central limit theorem (Theorem 2.4) is new.

**2. The model and the main results.** Throughout, let $p$ be a positive integer, let $b = b(p)$ and $n = n(p)$ be positive numbers depending on $p$, with $n$ an integer. (Many objects considered below depend on $p$, but we tend to suppress explicit reference to $p$ in the notation.) We assume the following concerning these numbers:

ASSUMPTION 2.1. As $p \to \infty$, we have $b \to \infty$, $n \to \infty$ and $b/n \to 0$, with $b \leq p$.

For any sequence of random variables $U_1, \ldots, U_n$, we let $\mathbf{C}(U_1, \ldots, U_n)$ denote their joint cumulant. (See Section 4 below for the definition of joint cumulants and a review of their properties.) Let

$$\{Z_j\}_{j=-\infty}^{\infty}$$

be a stationary sequence of real random variables, satisfying the following conditions:

ASSUMPTION 2.2.

(1) $$E(|Z_0|^k) < \infty \quad \text{for all } k \geq 1,$$

(2) $$EZ_0 = 0,$$

(3) $$\sum_{j_1} \cdots \sum_{j_r} |\mathbf{C}(Z_0, Z_{j_1}, \ldots, Z_{j_r})| < \infty \quad \text{for all } r \geq 1.$$

We refer to (3) as *joint cumulant summability*. In Section 2.4 below we describe a class of examples of sequences satisfying Assumption 2.2.

2.1. *Random matrices.* Let
$$\{\{Z_j^{(i)}\}_{j=-\infty}^{\infty}\}_{i=1}^{\infty}$$
be an i.i.d. family of copies of $\{Z_j\}_{j=-\infty}^{\infty}$. Let $X = X^{(p)}$ be the $n$-by-$p$ random matrix with entries
$$X(i,j) = X_{ij} = Z_j^{(i)}/\sqrt{n}.$$
Let $B = B^{(p)}$ be the $p$-by-$p$ deterministic matrix with entries
$$B(i,j) = B_{ij} = \begin{cases} 1, & \text{if } |i-j| \leq b, \\ 0, & \text{if } |i-j| > b. \end{cases}$$
Let $Y = Y^{(p)}$ be the $p$-by-$p$ random symmetric matrix with entries
$$Y(i,j) = Y_{ij} = B_{ij}(X^T X)_{ij} \tag{4}$$
and eigenvalues $\{\lambda_i^{(p)}\}_{i=1}^{p}$. Let
$$L = L^{(p)} = p^{-1} \sum_{i=1}^{p} \delta_{\lambda_i^{(p)}} \tag{5}$$
be the empirical measure of the eigenvalues of $Y$. Our attention will be focused on the limiting behavior of $L$ as $p \to \infty$.

2.2. *The measure $\nu_Z$.* For integers $j$ let
$$R(j) = \operatorname{Cov}(Z_0, Z_j). \tag{6}$$
Since $\mathbf{C}(Z_0, Z_j) = \operatorname{Cov}(Z_0, Z_j)$, a consequence of (3) is the existence of the *spectral density* $f_Z : [0,1] \to \mathbb{R}$ associated with the sequence $\{Z_j\}$, defined to be the Fourier transform
$$f_Z(\theta) = \sum_{j \in \mathbb{Z}} e^{2\pi i j \theta} R(j).$$
By the Szegö limit theorem [4], the empirical measure of the eigenvalues of the matrix $R(|i-j|)_{i,j=1}^{N}$ converges to the measure $\nu_Z := m \circ f_Z^{-1}$ on $\mathbb{R}$, where $m$ denotes Lebesgue measure on $[0,1]$. (Note that, considering the spectral density $f_Z$ as a random variable on the measure space $([0,1], m)$, one can interpret $\nu_Z$ as its law.) It is immediate to check from the definition that all moments of $\nu_Z$ are finite and are given by
$$\int_{\mathbb{R}} x^k \nu_Z(dx) = \int_0^1 f_Z(\theta)^k \, d\theta = \underbrace{R \star R \star \cdots \star R}_{k}(0)$$
$$= \sum_{\substack{i_1,\ldots,i_k \in \mathbb{Z} \\ i_1 + \cdots + i_k = 0}} \operatorname{Cov}(Z_0, Z_{i_1}) \cdots \operatorname{Cov}(Z_0, Z_{i_k}), \tag{7}$$

where $\star$ denotes convolution:

$$(F \star G)(j) = \sum_{k \in \mathbb{Z}} F(j-k)G(k),$$

for any two summable functions $F, G : \mathbb{Z} \to \mathbb{R}$. Note that (7) could just as well serve as the definition of $\nu_Z$.

2.3. *The coefficients $Q_{ij}$ and $R_i^{(m)}$.* With notation as in (3), (6), (7), for integers $m > 0$ and all integers $i$ and $j$, we write

(8)
$$Q_{ij} = \sum_{\ell \in \mathbb{Z}} \mathbf{C}(Z_i, Z_0, Z_{j+\ell}, Z_\ell),$$

$$R_i^{(m)} = \underbrace{R \star \cdots \star R}_{m}(i), \qquad R_i^{(0)} = \delta_{i0}.$$

By (3) the array $Q_{ij}$ is well defined and summable:

(9)
$$\sum_{i,j \in \mathbb{Z}} |Q_{ij}| < \infty.$$

The array $Q_{ij}$ is also symmetric:

(10)
$$Q_{ij} = \sum_{\ell \in \mathbb{Z}} \mathbf{C}(Z_{i-\ell}, Z_{-\ell}, Z_j, Z_0) = Q_{ji},$$

by stationarity of $\{Z_j\}$ and symmetry of $\mathbf{C}(\cdot, \cdot, \cdot, \cdot)$ under exchange of its arguments.

The following are the main results of this paper.

THEOREM 2.3 (Law of large numbers). *Let Assumptions 2.1 and 2.2 hold. Let $L = L^{(p)}$ be as in (5). Let $\nu_Z$ be as in (7). Then: $L$ converges weakly to $\nu_Z$, in probability.*

In other words, Theorem 2.3 implies that $L$ is a consistent estimator of $\nu_Z$, in the sense of weak convergence.

THEOREM 2.4 (Central limit theorem). *Let Assumptions 2.1 and 2.2 hold. Let $Y = Y^{(p)}$ be as in (4). Let $Q_{ij}$ and $R_i^{(m)}$ be as in (8). Then the process*

$$\left\{ \sqrt{\frac{n}{p}} (\operatorname{trace} Y^k - \mathbf{E} \operatorname{trace} Y^k) \right\}_{k=1}^{\infty}$$

*converges in distribution as $p \to \infty$ to a zero mean Gaussian process $\{G_k\}_{k=1}^{\infty}$ with covariance specified by the formula*

(11)
$$\frac{1}{k\ell} \mathbf{E} G_k G_\ell = 2R_0^{(k+\ell)} + \sum_{i,j \in \mathbb{Z}} R_i^{(k-1)} Q_{ij} R_j^{(\ell-1)}.$$

Note that the "correction" $Q_{ij}$ vanishes identically if $\{Z_j\}$ is Gaussian, compare Lemma 4.2 below.

2.4. *Some stationary sequences satisfying Assumption* 2.2. Fix a summable function $h:\mathbb{Z}\to\mathbb{R}$ and an i.i.d. sequence $\{W_\ell\}_{\ell=-\infty}^{\infty}$ of mean zero real random variables with moments of all orders. Now convolve: put $Z_j = \sum_\ell h(j+\ell)W_\ell$ for every $j$. It is immediate that (1) and (2) hold. To see the summability condition (3) on joint cumulants, assume at first that $h$ has finite support. Then, by standard properties of joint cumulants (the main point is covered by Lemma 4.1 below), we get the formula

$$(12) \qquad \mathbf{C}(Z_{j_0},\ldots,Z_{j_r}) = \sum_\ell h(j_0+\ell)\cdots h(j_r+\ell)\mathbf{C}(\underbrace{W_0,\ldots,W_0}_{r+1}),$$

which leads by a straightforward limit calculation to the analogous formula without the assumption of finite support of $h$, whence in turn verification of (3).

2.5. *Structure of the paper.* The proofs of Theorems 2.3 and 2.4 require a fair number of preliminaries. We provide them in the next few sections. In Section 3, we introduce some notation involving set partitions, and prove Proposition 3.1, which summarizes the properties of set partitions that we need. In spirit, if not in precise details, this section builds on [1]. In Section 4, we introduce joint cumulants and the Möbius inversion formula relating cumulants to moments, and in Section 5 we use the latter to calculate joint cumulants of random variables of the form trace $Y^k$ by manipulation of set partitions; see Proposition 5.2. In Section 6 we carry out some preliminary limit calculations in order to identify the dominant terms in the sums representing joint cumulants of random variables of the form trace $Y^k$. Finally, the proofs of Theorems 2.3 and 2.4 are completed in Sections 7 and 8, respectively.

## 3. A combinatorial estimate.

3.1. *Set partitions.* Given a positive integer $k$, we define $\mathrm{Part}(k)$ to be the family of subsets of the power set $2^{\{1,\ldots,k\}}$ consisting of sets $\Pi$ such that (i) $\varnothing \notin \Pi$, (ii) $\bigcup_{A\in\Pi} A = \{1,\ldots,k\}$, and (iii) for all $A,B\in\Pi$, if $A\neq B$, then $A\cap B = \varnothing$. Elements of $\mathrm{Part}(k)$ are called *set partitions* of $\{1,\ldots,k\}$, or context permitting simply *partitions*. Sometimes we call members of a partition *parts*. Given $\Pi, \Sigma \in \mathrm{Part}(k)$, we say that $\Sigma$ *refines* $\Pi$ (or is *finer* than $\Pi$) if for every $A\in\Sigma$ there exists some $B\in\Pi$ such that $A\subset B$. Given $\Pi, \Sigma \in \mathrm{Part}(k)$, let $\Pi\vee\Sigma \in \mathrm{Part}(k)$ be the least upper bound of $\Pi$ and $\Sigma$, that is, the finest partition refined by both $\Pi$ and $\Sigma$. We call $\Pi\in\mathrm{Part}(k)$ a *perfect matching* if every part of $\Pi$ has cardinality 2. Let $\text{Part}_2(k)$ be the subfamily of $\text{Part}(k)$ consisting of partitions $\Pi$ such that every part has cardinality at least 2. The cardinality of a set $S$ is denoted $\#S$, and $\lfloor x \rfloor$ denotes the greatest integer not exceeding $x$.

PROPOSITION 3.1.  *Let $k$ be a positive integer. Let $\Pi_0, \Pi_1, \Pi \in \text{Part}_2(2k)$ be given. Assume that $\Pi_0$ and $\Pi_1$ are perfect matchings. Assume that $\#\Pi_0 \vee \Pi_1 \vee \Pi = 1$. Then we have*

(13) $$\#\Pi_0 \vee \Pi + \#\Pi_1 \vee \Pi \leq 1 + \#\Pi \leq k + 1,$$

*and furthermore,*

(14) $$r > 1 \Rightarrow \#\Pi_0 \vee \Pi + \#\Pi_1 \vee \Pi \leq k + 1 - \lfloor r/2 \rfloor,$$

*where $r = \#\Pi_0 \vee \Pi_1$.*

The proposition is very close to [1], Lemma 4.10, almost a reformulation. But because the setup of [1] is rather different from the present one, the effort of translation is roughly equal to the effort of direct proof. We choose to give a direct proof in order to keep the paper self-contained. The proof will be finished in Section 3.5. In Section 9, we provide some comments concerning possible improvements of Proposition 3.1.

3.2. *Graphs.* We fix notation and terminology. The reader is encouraged to glance at Figure 1 when reading the rest of this section for an illustration of the various definitions in a concrete example.

3.2.1. *Basic definitions.* For us a *graph* $G = (V, E)$ is a pair consisting of a finite set $V$ and a subset $E \subset 2^V$ of the power set of $V$ such that every member of $E$ has cardinality 1 or 2. Elements of $V$ are called *vertices* and elements of $E$ are called *edges*. A *walk* $w$ on $G$ is a sequence $w = v_1 v_2 \cdots v_n$ of vertices of $G$ such that $\{v_i, v_{i+1}\} \in E$ for $i = 1, \ldots, n-1$, and in this situation we say that the initial point $v_1$ and terminal point $v_n$ of the walk are *joined* by $w$. A graph is *connected* if every two vertices are joined by a walk. For any connected graph, $\#V \leq 1 + \#E$. A graph $G = (V, E)$ is called a *tree* if connected and further $\#V = 1 + \#E$. Alternatively, a connected graph $G = (V, E)$ is a tree if and only if there exists no edge $e \in E$ such that the subgraph $G' = (V, E \setminus \{e\})$ gotten by "erasing" the edge $e$ is connected.

For future reference, we quote without proof the following elementary lemma.

LEMMA 3.2 (Parity principle).  *Let $w = v_1 \cdots v_n$ be a walk on a tree $T = (V, E)$ beginning and ending at the same vertex, that is, such $v_1 = v_n$. Then $w$ visits every edge of $T$ an even number of times, that is,*

$$\#\{i \in \{1, \ldots, n-1\} \mid \{v_i, v_{i+1}\} = e\}$$

*is an even number for every $e \in E$.*

3.3. *Reduction of $\Pi_0$ and $\Pi_1$ to standard form.* After relabeling the elements of $\{1, \ldots, 2k\}$, we may assume that for some positive integers $k_1, \ldots, k_r$ summing to $k$ we have

$$\Pi_0 \vee \Pi_1 = \{(K_{\alpha-1}, K_\alpha] \cap \mathbb{Z} \mid \alpha = 1, \ldots, r\},$$

where $K_\alpha = 2\sum_{\beta < \alpha} k_\beta$ for $\alpha = 0, \ldots, r$, and after some further relabeling, we may assume that

$$\Pi_0 = \{\{2i - 1, 2i\} \mid i = 1, \ldots, k\}.$$

It is well known (and easily checked) that for any perfect matchings $\Sigma_0, \Sigma_1 \in \mathrm{Part}_2(2k)$, the graph $(\{1, \ldots, 2k\}, \Sigma_0 \cup \Sigma_1)$ is a disjoint union of $\#\Sigma_0 \vee \Sigma_1$ graphs of the form

$$(\{1, 2\}, \{\{1, 2\}\}), \qquad (\{1, 2, 3, 4\}, \{\{1, 2\}, \{2, 3\}, \{3, 4\}, \{4, 1\}\}),$$
$$(\{1, 2, 3, 4, 5, 6\}, \{\{1, 2\}, \{2, 3\}, \{3, 4\}, \{4, 5\}, \{5, 6\}, \{6, 1\}\})$$

and so on. (The intuition is that the members of $\Sigma_0$ and $\Sigma_1$ "join hands" alternately to form cycles.) Thus, after a final round of relabeling, we may assume that

$$\Pi_1 = \bigcup_{\alpha=1}^r (\{\{i_{2k_\alpha}^{(\alpha)}, i_1^{(\alpha)}\}\} \cup \{\{i_{2\nu}^{(\alpha)}, i_{2\nu+1}^{(\alpha)}\} \mid \nu = 1, \ldots, k_\alpha - 1\}),$$

where $i_\nu^{(\alpha)} = K_{\alpha-1} + \nu$. Note that

$$\Pi_0 = \bigcup_{\alpha=1}^r \{\{i_{2\nu-1}^{(\alpha)}, i_{2\nu}^{(\alpha)}\} \mid \nu = 1, \ldots, k_\alpha\},$$

$$\Pi_0 \vee \Pi_1 = \{\{i_1^{(\alpha)}, \ldots, i_{2k_\alpha}^{(\alpha)}\} \mid \alpha = 1, \ldots, r\}$$

in terms of the notation introduced to describe $\Pi_1$.

3.4. *Graph-theoretical "coding" of $\Pi$.*

3.4.1. *Construction of a graph $G$.* For $i = 0, 1$, let

$$\varphi_i \colon \{1, \ldots, 2k\} \to V_i$$

be an onto function such that $\Pi_i \vee \Pi$ is the family of level sets for $\varphi_i$. Assume further that $V_0 \cap V_1 = \varnothing$. We now define a graph $G = (V, E)$ by declaring that

$$V = V_0 \cup V_1,$$
$$E = \{\{\varphi_0(i), \varphi_1(i)\} \mid i = 1, \ldots, 2k\}.$$

LEMMA 3.3. *$G$ is connected.*

Because $\varphi_i(j) = \varphi_i(\ell)$ for $i = 0,1$ if $j, \ell$ belong to the same part of $\Pi$, we must have $\#E \leq \#\Pi$. Further, $\#\Pi \leq k$ since $\Pi \in \text{Part}_2(2k)$. Thus, using Lemma 3.3 in the first inequality, we have

$$\#\Pi_0 \vee \Pi + \#\Pi_1 \vee \Pi = \#V \leq 1 + \#E \leq 1 + \#\Pi \leq k+1,$$

which proves inequality (13) of Proposition 3.1.

PROOF OF LEMMA 3.3. Suppose rather that we have a decomposition $V = X \cup Y$ where $X \cap Y = \varnothing$, $X \neq \varnothing$, $Y \neq \varnothing$, and no edge of $G$ joins a vertex in $X$ to a vertex in $Y$. Consider the subsets

$$I = \varphi_0^{-1}(V_0 \cap X) \cup \varphi_1^{-1}(V_1 \cap X), \qquad J = \varphi_0^{-1}(V_0 \cap Y) \cup \varphi_1^{-1}(V_1 \cap Y)$$

of $\{1, \ldots, 2k\}$. Clearly $I \cup J = \{1, \ldots, 2k\}$, $I \neq \varnothing$, and $J \neq \varnothing$. We claim that $I \cap J = \varnothing$. Suppose rather that there exists $i \in I \cap J$. Then we must either have $\varphi_0(i) \in V_0 \cap X$ and $\varphi_1(i) \in V_1 \cap Y$, or else $\varphi_1(i) \in V_1 \cap X$ and $\varphi_0(i) \in V_0 \cap Y$. In either case we have exhibited an edge of $G$ connecting a vertex in $X$ to a vertex in $Y$, which is a contradiction. Therefore $I \cap J = \varnothing$. Thus the set $\{I, J\} \in \text{Part}(2k)$ is a partition refined by both $\Pi_0 \vee \Pi$ and $\Pi_1 \vee \Pi$, which is a contradiction to $\#\Pi_0 \vee \Pi_1 \vee \Pi = 1$. Therefore $G$ is connected. □

LEMMA 3.4. *There exist walks*
$$w^{(\alpha)} = v_1^{(\alpha)} \cdots v_{2k_\alpha+1}^{(\alpha)} \qquad \text{for } \alpha = 1, \ldots, r$$
*on $G$ such that*
$$v_1^{(\alpha)} = v_{2k_\alpha+1}^{(\alpha)},$$
$$\{\varphi_0(i_\nu^{(\alpha)}), \varphi_1(i_\nu^{(\alpha)})\} = \begin{cases} \{v_\nu^{(\alpha)}, v_{\nu+1}^{(\alpha)}\}, & \text{if } \nu < 2k_\alpha, \\ \{v_{2k_\alpha}^{(\alpha)}, v_1^{(\alpha)}\}, & \text{if } \nu = 2k_\alpha, \end{cases}$$
*for $\alpha = 1, \ldots, r$ and $\nu = 1, \ldots, 2k_\alpha$.*

PROOF. We define
$$v_\nu^{(\alpha)} = \begin{cases} \varphi_1(i_\nu^{(\alpha)}), & \text{if } \nu \text{ is odd and } \nu < 2k_\alpha, \\ \varphi_0(i_\nu^{(\alpha)}), & \text{if } \nu \text{ is even}, \\ \varphi_1(i_1^{(\alpha)}), & \text{if } \nu = 2k_\alpha + 1, \end{cases}$$

for $\alpha = 1, \ldots, r$ and $\nu = 1, \ldots, 2k_\alpha + 1$. Clearly we have $v_1^{(\alpha)} = v_{2k_\alpha+1}^{(\alpha)}$. Recalling that $\varphi_0$ by construction is constant on the set $\{i_1^{(\alpha)}, i_2^{(\alpha)}\} \in \Pi_0$, we see that

$$\{\varphi_1(i_1^{(\alpha)}), \varphi_0(i_1^{(\alpha)})\} = \{\varphi_1(i_1^{(\alpha)}), \varphi_0(i_2^{(\alpha)})\} = \{v_1^{(\alpha)}, v_2^{(\alpha)}\}.$$

By similar considerations one checks the remaining claims of the lemma. We omit further details. □

LEMMA 3.5. *Assume that $r > 1$. For every $A \in \Pi_0 \vee \Pi_1$ there exists an index $m \in A$, a set $A' \in \Pi_0 \vee \Pi_1$ distinct from $A$ and an index $m' \in A'$ such that $\{\varphi_0(m), \varphi_1(m)\} = \{\varphi_0(m'), \varphi_1(m')\}$.*

In other words, if $r > 1$, then for every walk $w^{(\alpha)}$, there is an edge $e$ of $G$ and another walk $w^{(\alpha')}$ such that both $w^{(\alpha)}$ and $w^{(\alpha')}$ visit $e$.

PROOF OF LEMMA 3.5. Because $\#\Pi_0 \vee \Pi_1 \vee \Pi = 1$, given $A \in \Pi_0 \vee \Pi_1$, there must exist $A' \in \Pi_0 \vee \Pi_1$ distinct from $A$ and a set $B \in \Pi$ such that $A \cap B \neq \varnothing$ and $A' \cap B \neq \varnothing$. Choose $m \in A \cap B$ and $m' \in A' \cap B$. Because the functions $\varphi_0$ and $\varphi_1$ are constant on the set $B$, we are done. $\square$

3.5. *Completion of the proof of Proposition* 3.1. We have seen that Lemma 3.3 proves inequality (13). We just have to prove inequality (14). Assume that $r > 1$ for the rest of the proof. Consider the graph $G = (V, E)$ as in Section 3.4. Let $E' \subset E$ be such that $T = (V, E')$ is a tree (such a choice is possible because $G$ is connected). It will be enough to show that $\#E' \leq k - r/2$. Now we adapt to the present situation a device ("edge-bounding tables") introduced in the proof of [1], Lemma 4.10. Let us call a function $f : \{1, \ldots, 2k\} \to \{0, 1\}$ a *good estimator* under the following conditions:

- For all $i \in \{1, \ldots, 2k\}$, if $f(i) = 1$, then $\{\varphi_0(i), \varphi_1(i)\} \in E'$.
- For each $e \in E'$ there exist distinct $i, j \in \{1, \ldots, 2k\}$ such that $e = \{\varphi_0(i), \varphi_1(i)\} = \{\varphi_0(j), \varphi_1(j)\}$ and $f(i) = f(j) = 1$.
- For each $e \in E'$ and $A \in \Pi_0 \vee \Pi_1$, if there exists $\ell \in A$ such that $e = \{\varphi_0(\ell), \varphi_1(\ell)\}$, then there exists $\ell' \in A$ such that $e = \{\varphi_0(\ell'), \varphi_1(\ell')\}$ and $f(\ell') = 1$.

For a good estimator $f$ we automatically have $\frac{1}{2} \sum f(i) \geq \#E'$. By definition a good estimator is bounded above by the indicator of the set $\{i \in \{1, \ldots, 2k\} \mid \{\varphi_0(i), \varphi_1(i)\} \in E'\}$, and such an indicator function is an example of a good estimator. Fix now any good estimator $f$. Suppose that on some set $A = \{i_1^{(\alpha)}, \ldots, i_{2k_\alpha}^{(\alpha)}\} \in \Pi_0 \vee \Pi_1$ the function $f$ is identically equal to 1. Then the corresponding walk $w^{(\alpha)}$ on $G$ is a walk on $T$, and by the Parity Principle (Lemma 3.2) visits every edge of $T$ an even number of times. Select $m \in A$ as in Lemma 3.5. Let $g$ be the function agreeing with $f$ everywhere except that $g(m) = 0$. Then $g$ is again a good estimator. Continuing in this way we can construct a good estimator not identically equal to 1 on any of the sets $A \in \Pi_0 \vee \Pi_1$, whence the desired estimate $\#E \leq k - r/2$.

Figure 1 illustrates the various objects studied in this section.

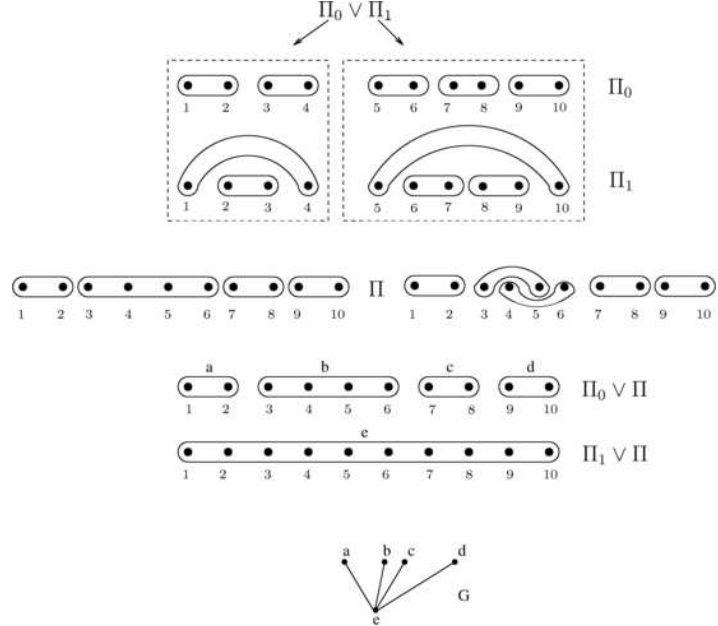

FIG. 1. *Two different partitions $\Pi$ for which $k = 5$, $k_1 = 2$, $k_2 = 3$, such that both are associated to the same graph $G = (V, E)$, where $V = \{a, b, c, d, e\}$. Note that both partitions generate walks eaebe and ebecede on $G$.*

## 4. Joint cumulants.

4.1. *Definition.* Let $X_1, \ldots, X_k$ be real random variables defined on a common probability space with moments of all orders, in which case the characteristic function $\mathbf{E} \exp(\sum_{j=1}^{k} it_j X_j)$ is an infinitely differentiable function of the real variables $t_1, \ldots, t_k$. One defines the *joint cumulant* $\mathbf{C}(X_1, \ldots, X_k)$ by the formula

$$\mathbf{C}(X_1, \ldots, X_k) = \mathbf{C}\{X_i\}_{i=1}^{k}$$
$$= i^{-k} \frac{\partial^k}{\partial t_1 \cdots \partial t_k} \log \mathbf{E} \exp\left(\sum_{j=1}^{k} it_j X_j\right)\bigg|_{t_1 = \cdots = t_k = 0}.$$

(The middle expression is a convenient abbreviated notation.) The quantity $\mathbf{C}(X_1, \ldots, X_k)$ depends symmetrically and $\mathbb{R}$-multilinearly on $X_1, \ldots, X_k$. Moreover, dependence is continuous with respect to the $L^k$-norm. One has in particular

$$\mathbf{C}(X) = \mathbf{E} X, \qquad \mathbf{C}(X, X) = \operatorname{Var} X, \qquad \mathbf{C}(X, Y) = \operatorname{Cov}(X, Y).$$

The following standard properties of joint cumulants will be used. Proofs are omitted.

LEMMA 4.1. *If there exists $0 < \ell < k$ such that the $\sigma$-fields $\sigma\{X_i\}_{i=1}^{\ell}$ and $\sigma\{X_i\}_{i=\ell+1}^{k}$ are independent, then $\mathbf{C}(X_1,\ldots,X_k) = 0$.*

LEMMA 4.2. *The random vector $X_1,\ldots,X_k$ has a Gaussian joint distribution if and only if $\mathbf{C}(X_{i_1},\ldots,X_{i_r}) = 0$ for every integer $r \geq 3$ and sequence $i_1,\ldots,i_r \in \{1,\ldots,k\}$.*

4.2. *Combinatorial description of joint cumulants.* As above, let $X_1,\ldots,X_k$ be real random variables defined on a common probability space with moments of all orders. Let $\Pi \in \mathrm{Part}(k)$ also be given. We define

$$\mathbf{C}_\Pi(X_1,\ldots,X_k) = \mathbf{C}_\Pi\{X_i\}_{i=1}^k = \prod_{A \in \Pi} \mathbf{C}\{X_i\}_{i \in A},$$

$$\mathbf{E}_\Pi(X_1,\ldots,X_k) = \mathbf{E}_\Pi\{X_i\}_{i=1}^k = \prod_{A \in \Pi} \mathbf{E} \prod_{i \in A} X_i.$$

(The middle expressions are convenient abbreviations.) Note that if $X_1,\ldots,X_k$ are zero mean random variables, then $\mathbf{C}_\Pi(X_1,\ldots,X_k)$ vanishes unless $\Pi \in \mathrm{Part}_2(k)$. The formula

$$\mathbf{E} X_1 \cdots X_k = \sum_{\Pi \in \mathrm{Part}(k)} \mathbf{C}_\Pi(X_1,\ldots,X_k) \tag{15}$$

is well-known, and can be verified in a straightforward way by manipulating Taylor expansions of characteristic functions. More generally we have the following lemma, whose proof can be found in [6], page 290.

LEMMA 4.3. *With $X_1,\ldots,X_k$ as above, and for all $\Pi \in \mathrm{Part}(k)$, we have*

$$\mathbf{E}_\Pi\{X_i\}_{i=1}^k = \sum_{\substack{\Sigma \in \mathrm{Part}(k) \\ \Sigma \text{ refines } \Pi}} \mathbf{C}_\Sigma\{X_i\}_{i=1}^k, \tag{16}$$

$$\mathbf{C}_\Pi\{X_i\}_{i=1}^k$$
$$= \sum_{\substack{\Sigma \in \mathrm{Part}(k) \\ \Sigma \text{ refines } \Pi}} \left( \prod_{A \in \Pi} (-1)^{\#\{B \in \Sigma | B \subset A\} - 1}(\#\{B \in \Sigma \mid B \subset A\} - 1)! \right) \tag{17}$$
$$\times \mathbf{E}_\Sigma\{X_i\}_{i=1}^k.$$

We will use the following algebraic fact to compute joint cumulants. For a proof see, for example, [7], Example 3.10.4.

LEMMA 4.4 [Möbius Inversion for the poset $\operatorname{Part}(k)$]. *Let $A$ be an Abelian group and let $f, g : \operatorname{Part}(k) \to A$ be functions. Then we have*

$$\tag{18} \left( [\forall \Sigma \in \operatorname{Part}(k)] \ f(\Sigma) = \sum_{\substack{\Pi \in \operatorname{Part}(k) \\ \Pi \text{ refines } \Sigma}} g(\Pi) \right)$$

*if and only if*

$$\tag{19} \left( [\forall \Pi \in \operatorname{Part}(k)] \quad g(\Pi) = \sum_{\substack{\Sigma \in \operatorname{Part}(k) \\ \Sigma \text{ refines } \Pi}} \left( \prod_{A \in \Pi} (-1)^{\#\{B \in \Sigma | B \subset A\} - 1} (\#\{B \in \Sigma \mid B \subset A\} - 1)! \right) f(\Sigma) \right).$$

In applications below we will simply have $A = \mathbb{R}$.

**5. Cumulant calculations.** In the context of matrix models, cumulants are useful because they allow one to replace enumeration over arbitrary graphs by enumeration over *connected* graphs. We wish to mimic this idea in our context. We first describe the setup, and then perform some computations that culminate in Proposition 5.2, which gives an explicit formula for joint cumulants of random variables of the form $\operatorname{trace} Y^k$.

5.1. *The setup.* An $(n, k)$-*word* $\mathbf{i}$ is by definition a function

$$\mathbf{i} : \{1, \ldots, k\} \to \{1, \ldots, n\}.$$

Given $\Pi \in \operatorname{Part}(k)$ and an $(n, k)$-word $\mathbf{i}$, we say that $\mathbf{i}$ is $\Pi$-*measurable* if $\mathbf{i}$ is constant on each set belonging to $\Pi$. Similarly and more generally, we speak of the $\Pi$-measurability of any function $\mathbf{i} : \{1, \ldots, k\} \to \mathbb{Z}$.

Let $r$ be a positive integer. Let $k_1, \ldots, k_r$ be positive integers and put $k = k_1 + \cdots + k_r$. Let special perfect matchings $\Pi_0, \Pi_1 \in \operatorname{Part}(2k)$ be defined as follows:

$$\Pi_0 = \{\{1, 2\}, \{3, 4\}, \ldots, \{2k-3, 2k-2\}, \{2k-1, 2k\}\},$$
$$\Pi_1 = \{\{2, 3\}, \ldots, \{K_1, 1\}, \{K_1 + 2, K_1 + 3\}, \ldots, \{K_2, K_1 + 1\},$$
$$\ldots, \{K_{r-1} + 2, K_{r-1} + 3\}, \ldots, \{K_r, K_{r-1} + 1\}\},$$

where $K_i = 2\sum_{j=1}^i k_j$ for $i=1,\ldots,r$. (Thus $\Pi_0$ and $\Pi_1$ are in the standard form discussed in Section 3.3 above.) To abbreviate, for any $\Pi \in \text{Part}(2k)$ and $(p, 2k)$-word $\mathbf{j}$, put

$$B(\mathbf{j}) = \prod_{\alpha=1}^{2k} B(\mathbf{j}(2\alpha-1), \mathbf{j}(2\alpha)), \qquad C_\Pi(\mathbf{j}) = \mathbf{C}_\Pi(Z_{\mathbf{j}(1)}, \ldots, Z_{\mathbf{j}(2k)}).$$

Note that, on the one hand, $B(\mathbf{j})$ depends on $p$ even though the notation does not show the dependence. Note that, on the other hand, $C_\Pi(\mathbf{j})$ is independent of $p$. Indeed, $C_\Pi(\mathbf{j})$ remains well defined by the formula above for any function $\mathbf{j}: \{1,\ldots,2k\} \to \mathbb{Z}$.

Concerning the numbers $C_\Pi(\mathbf{j})$ we record for later reference the following consequence of the joint cumulant summability hypothesis (3) and the stationarity of $\{Z_j\}$. The proof is immediate from the definitions and therefore omitted. Let $\mathbb{Z}^\Pi$ be the subgroup of $\mathbb{Z}^{2k}$ consisting of functions on $\{1,\ldots,2k\}$ constant on each part of $\Pi$.

LEMMA 5.1. *For every* $\mathbf{j}: \{1,\ldots,2k\} \to \mathbb{Z}$, (i) *the value of* $C_\Pi(\mathbf{j})$ *depends only on the coset of* $\mathbb{Z}^\Pi$ *to which* $\mathbf{j}$ *belongs and moreover* (ii) *we have*

$$\sum_{\mathbf{j} \in \mathbb{Z}^{2k}/\mathbb{Z}^\Pi} |C_\Pi(\mathbf{j})| < \infty.$$

The lemma will be the basis for our limit calculations.
Our immediate goal is to prove the following result.

PROPOSITION 5.2. *With the previous notation, we have*

$$\mathbf{C}(\text{trace}\, Y^{k_1}, \ldots, \text{trace}\, Y^{k_r})$$

(20)
$$= \sum_{\substack{\Pi \in \text{Part}_2(2k) \\ \text{s.t.}\ \#\Pi_0 \vee \Pi_1 \vee \Pi = 1}} n^{-k + \#\Pi_0 \vee \Pi} \sum_{\substack{\mathbf{j}:\, (p,2k)\text{-word s.t.}\ \mathbf{j} \\ \text{is}\ \Pi_1\text{-measurable}}} B(\mathbf{j}) C_\Pi(\mathbf{j}).$$

PROOF. The proof involves an application of the Möbius inversion formula (Lemma 4.4). Recall that

$$\text{trace}\, Y^k = \sum_{\mathbf{l}:\, (n,k)\text{-word}} Y(\mathbf{l}(1), \mathbf{l}(2)) Y(\mathbf{l}(2), \mathbf{l}(3)) \cdots Y(\mathbf{l}(k), \mathbf{l}(1)).$$

Further, $Y(j_1, j_2) = B(j_1, j_2) \sum_i X(i, j_1) X(i, j_2)$ and hence

$$\text{trace}\, Y^k = \sum_{\substack{\mathbf{i}:\, (n,2k)\text{-word s.t.} \\ \mathbf{i}(2t-1)=\mathbf{i}(2t),\, t=1,\ldots,k,\ \text{and} \\ \mathbf{j}:\, (n,2k)\text{-word s.t.} \\ \mathbf{j}(2t)=\mathbf{j}(2t+1),\, t=1,\ldots,k}} B(\mathbf{j}) \prod_{\alpha=1}^{2k} X(\mathbf{i}(\alpha), \mathbf{j}(\alpha)),$$

where $\mathbf{j}(2k+1)$ is defined by the "wrap-around rule": $\mathbf{j}(2k+1) = \mathbf{j}(1)$. Hence, we have

$$\mathbf{E}(\operatorname{trace} Y^{k_1}) \cdots (\operatorname{trace} Y^{k_r})$$

(21)
$$= \sum_{\substack{\mathbf{i}:\,(n,2k)\text{-word s.t. } \mathbf{ij}:\,(p,2k)\text{-word s.t. } \mathbf{j} \\ \text{is } \Pi_0\text{-measurable} \quad \text{is}\Pi_1\text{-measurable}}} B(\mathbf{j}) \mathbf{E} \prod_{\alpha=1}^{2k} X(\mathbf{i}(\alpha), \mathbf{j}(\alpha))$$

Using the representation of moments in terms of cumulants [see (15)] we get

$$\mathbf{E}(\operatorname{trace} Y^{k_1}) \cdots (\operatorname{trace} Y^{k_r})$$

$$= \sum_{\substack{\mathbf{i}:\,(n,2k)\text{-word s.t. } \mathbf{ij}:\,(p,2k)\text{-word s.t. } \mathbf{j} \\ \text{is } \Pi_0\text{-measurable} \quad \text{is } \Pi_1\text{-measurable}}} B(\mathbf{j})$$

$$\times \sum_{\Pi \in \operatorname{Part}(2k)} \mathbf{C}_\Pi \{X(\mathbf{i}(\alpha), \mathbf{j}(\alpha))\}_{\alpha=1}^{2k}$$

(22)
$$= \sum_{\substack{\mathbf{i}:\,(n,2k)\text{-word s.t. } \mathbf{ij}:\,(p,2k)\text{-word s.t. } \mathbf{j} \\ \text{is } \Pi_0\text{-measurable} \quad \text{is } \Pi_1\text{-measurable}}} n^{-k} B(\mathbf{j}) \sum_{\substack{\Pi \in \operatorname{Part}_2(2k) \\ \text{s.t. } \mathbf{i} \text{ is } \Pi\text{-measurable}}} C_\Pi(\mathbf{j})$$

$$= \sum_{\Pi \in \operatorname{Part}_2(2k)} n^{-k + \#\Pi_0 \vee \Pi} \sum_{\substack{\mathbf{j}:\,(p,2k)\text{-word} \\ \text{s.t. } \mathbf{j} \text{ is } \Pi_1\text{-measurable}}} B(\mathbf{j}) C_\Pi(\mathbf{j}),$$

where in the next to last equality we used that cumulants of independent variables vanish in order to restrict the summation to words $\mathbf{i}$ that are $\Pi$-measurable.

We next define an embedding of $\operatorname{Part}(r)$ in $\operatorname{Part}(2k)$. It will be convenient to use $\pi, \sigma$ to denote elements of $\operatorname{Part}(r)$ and $\Pi, \Sigma$ to denote elements of $\operatorname{Part}(2k)$. (Also we use upper case Roman letters for subsets of $\{1, \ldots, 2k\}$ and lower case Roman letters for subsets of $\{1, \ldots, r\}$.) Put

$$A_1 = \{1, \ldots, K_1\}, \ldots, A_r = \{K_{r-1} + 1, \ldots, K_r\},$$

so that

$$\Pi_0 \vee \Pi_1 = \{A_1, \ldots, A_r\}.$$

Given $a \subset \{1, \ldots, r\}$, let $a^* = \bigcup_{i \in a} A_i$, and given $\sigma \in \operatorname{Part}(r)$, let

$$T(\sigma) = \{a^* \mid a \in \sigma\} \in \operatorname{Part}(2k).$$

Via $T$ the poset $\operatorname{Part}(r)$ maps isomorphically to the subposet of $\operatorname{Part}(2k)$ consisting of partitions refined by $\Pi_0 \vee \Pi_1$.

We are ready to apply the Möbius inversion formula (Lemma 4.4). Consider the real-valued functions $f$ and $g$ on $\mathrm{Part}(r)$ defined as follows:

$$(23) \quad g(\pi) = \sum_{\substack{\Pi \in \mathrm{Part}(2k) \\ \Pi_0 \vee \Pi_1 \vee \Pi = T(\pi)}} n^{-k + \#\Pi_0 \vee \Pi} \sum_{\substack{\mathbf{j} \colon (p, 2k)\text{-word} \\ \text{s.t. } \mathbf{j} \text{ is } \Pi_1\text{-measurable}}} B(\mathbf{j}) C_\Pi(\mathbf{j})$$

and

$$(24) \quad f(\sigma) = \sum_{\substack{\pi \in \mathrm{Part}(r) \\ \pi \text{ refines } \sigma}} g(\pi).$$

Now $\pi$ refines $\sigma$ if and only if $T(\pi)$ refines $T(\sigma)$. Therefore we have

$$(25) \quad f(\sigma) = \sum_{\substack{\Pi \in \mathrm{Part}(2k) \\ \Pi_0 \vee \Pi_1 \vee \Pi \text{ refines } T(\sigma)}} n^{-k + \#\Pi_0 \vee \Pi} \sum_{\substack{\mathbf{j} \colon (p, 2k)\text{-word} \\ \text{s.t. } \mathbf{j} \text{ is } \Pi_1\text{-measurable}}} B(\mathbf{j}) C_\Pi(\mathbf{j}).$$

Using (24) and applying Lemma 4.4, it follows that for any $\pi \in \mathrm{Part}(r)$,

$$(26) \quad g(\pi) = \sum_{\substack{\sigma \in \mathrm{Part}(r) \\ \sigma \text{ refines } \pi}} \left( \prod_{a \in \pi} (-1)^{\#\{b \in \sigma | b \subset a\} - 1} (\#\{b \in \sigma \mid b \subset a\} - 1)! \right) f(\sigma).$$

An evident modification of the calculation (22) above gives for every $\sigma \in \mathrm{Part}(r)$ that $\mathbf{E}_\sigma(\mathrm{trace}\, Y^{k_1}, \ldots, \mathrm{trace}\, Y^{k_r})$ equals the right-hand side of (25), and therefore equals $f(\sigma)$. Thus, (26), when compared with (17), shows that

$$g(\{\{1, \ldots, r\}\}) = \mathbf{C}(\mathrm{trace}\, Y^{k_1}, \ldots, \mathrm{trace}\, Y^{k_r}),$$

which is exactly what we wanted to prove. $\square$

**6. Limit calculations.** We continue in the setting of Proposition 5.2. We find the order of magnitude of the subsum of the right-hand side of (20) indexed by $\Pi$ and compute limits as $p \to \infty$ in certain cases.

PROPOSITION 6.1. *Fix $\Pi \in \mathrm{Part}_2(2k)$ such that $\#\Pi_0 \vee \Pi_1 \vee \Pi = 1$. We have*

$$(27) \quad \sum_{\substack{\mathbf{j} \colon (p, 2k)\text{-word s.t. } \mathbf{j} \\ \text{is } \Pi_1\text{-measurable}}} B(\mathbf{j}) C_\Pi(\mathbf{j}) = O_{p \to \infty}(p b^{-1 + \#\Pi_1 \vee \Pi})$$

*where the implied constant depends only on $\Pi_0$, $\Pi_1$ and $\Pi$.*

Before commencing the proof of the proposition, we record an elementary lemma which expresses in algebraic terms the fact that a tree is connected and simply connected. We omit the proof. We remark that a tree can have no edges joining a vertex to itself.

LEMMA 6.2. *Let $T = (V, E)$ be a tree with vertex set $V \subset \{1, \ldots, 2k\}$. For each function $\mathbf{j}: V \to \mathbb{Z}$ define $\delta\mathbf{j}: E \to \mathbb{Z}$ by the rule*

$$\delta\mathbf{j}(\{\alpha, \beta\}) = \mathbf{j}(\beta) - \mathbf{j}(\alpha)$$

*for all $\alpha, \beta \in V$ such that $\alpha < \beta$ and $\{\alpha, \beta\} \in E$. Then:* (i) $\delta\mathbf{j} = 0$ *implies that $\mathbf{j}$ is constant.* (ii) *For every $\mathbf{k}: E \to \mathbb{Z}$ there exists $\mathbf{j}: V \to \mathbb{Z}$ unique up to addition of a constant such that $\delta\mathbf{j} = \mathbf{k}$.*

We will refer to $\delta$ as the *increment operator* associated to the tree $T$.

PROOF OF PROPOSITION 6.1. We begin by constructing a tree $T$ to which Lemma 6.2 will be applied. Let $\tilde{E}_2$ be the set consisting of all two-element subsets of parts of $\Pi$. With

$$V = \{1, \ldots, 2k\},$$

consider the graphs

$$G_{012} = (V, \Pi_0 \cup \Pi_1 \cup \tilde{E}_2), \qquad G_{12} = (V, \Pi_1 \cup \tilde{E}_2), \qquad G_2 = (V, \tilde{E}_2).$$

By hypothesis the graph $G_{012}$ is connected, and further, the number of connected components of $G_{12}$ (resp., $G_2$) equals $\#\Pi_1 \vee \Pi$ (resp., $\#\Pi$). Now choose $E_2 \subset \tilde{E}_2$ so that $T_2 = (V, E_2)$ is a spanning forest in $G_2$, that is, a subgraph with the same vertices but the smallest number of edges possible consistent with having the same number of connected components. Then choose $E_1 \subset \Pi_1$ such that $T_{12} = (V, E_1 \cup E_2)$ is a spanning forest in $G_{12}$, and finally choose $E_0 \subset \Pi_0$ such that $T_{012} = (V, E_0 \cup E_1 \cup E_2)$ is a spanning tree in $G_{012}$. By construction, the sets $E_i$, $i = 0, 1, 2$, are disjoint. Note that Lemma 6.2 applies not only to $T_{012}$, but also to the connected components of $T_{12}$ and $T_2$. Note that

(28) $$\#E_0 = -1 + \#\Pi_1 \vee \Pi$$

by construction. Hereafter we write simply $T = T_{012}$.

The bound in (27) will be obtained by relaxing some of the constraints concerning the collection of words $\mathbf{j}$ over which the summation runs. We will work with the increment operator $\delta$ associated to $T$ by Lemma 6.2. For $i = 0, 1, 2$ let $S_i$ be the Abelian group (independent of $p$) consisting of functions $\mathbf{j}: V \to \mathbb{Z}$ such that:

- $\mathbf{j}(1) = 0$,
- $\delta\mathbf{j}$ is supported on the set $E_i$.

Also let

$$S_{-1} = \{\mathbf{j}: V \to \mathbb{Z} \mid \delta\mathbf{j} = 0\} = \{\mathbf{j}: V \to \mathbb{Z} \mid \mathbf{j}: \text{constant}\},$$

which is independent of $p$. Recall that for any partition $\Pi$, $\mathbb{Z}^\Pi$ is the subgroup of $\mathbb{Z}^{2k}$ consisting of functions on $\{1,\ldots,2k\}$ constant on each part of $\Pi$. By Lemma 6.2 applied to $T$ and also to the connected components of $T_{12}$ and $T_2$, we have

$$
\begin{aligned}
\mathbb{Z}^{2k} &= S_{-1} \oplus S_0 \oplus S_1 \oplus S_2, \\
\mathbb{Z}^\Pi &= S_{-1} \oplus S_0 \oplus S_1, \\
\mathbb{Z}^{\Pi_1 \vee \Pi} &= S_{-1} \oplus S_0.
\end{aligned}
\tag{29}
$$

Let $S_0^{(p)} \subset S_{-1} \oplus S_0$ be the subset (depending on $p$) consisting of functions $\mathbf{j}: V \to \mathbb{Z}$ such that:

- $\mathbf{j}(1) \in \{1,\ldots,p\}$,
- $|\delta\mathbf{j}(e)| \leq b$ for all $e \in E_0$.

Now if $\mathbf{j}$ is a $\Pi_1$-measurable $(p,2k)$-word such that $B(\mathbf{j})$ does not vanish, then the following hold:

(30)
- $\mathbf{j}(1) \in \{1,\ldots,p\}$,
- $|\delta\mathbf{j}(e)| \leq b$ for $e \in E_0$ (because $E_0 \subset \Pi_0$),
- $\delta\mathbf{j}(e) = 0$ for $e \in E_1$ (because $E_1 \subset \Pi_1$).

By (29) it follows that a $\Pi_1$-measurable $(p,2k)$-word $\mathbf{j}$ such that $B(\mathbf{j})$ does not vanish has a unique decomposition $\mathbf{j} = \mathbf{j}_0 + \mathbf{j}_2$ with $\mathbf{j}_0 \in S_0^{(p)}$ and $\mathbf{j}_2 \in S_2$, and moreover we necessarily have

$$C_\Pi(\mathbf{j}) = C_\Pi(\mathbf{j}_2) \tag{31}$$

by Lemma 5.1(i) and the $\Pi$-measurability of $\mathbf{j}_0$.

We now come to the end of the proof. We have

$$\sum_{\substack{\mathbf{j}:\, (p,2k)\text{-word} \\ \text{s.t. } \mathbf{j} \text{ is} \\ \Pi_1\text{-measurable}}} |B(\mathbf{j})C_\Pi(\mathbf{j})| \leq \#S_0^{(p)} \sum_{\mathbf{j} \in S_2} |C_\Pi(\mathbf{j})|$$

$$\leq p(2b+1)^{-1+\#\Pi_1 \vee \Pi} \sum_{\mathbf{j} \in S_2} |C_\Pi(\mathbf{j})|$$

at the first inequality by (29), (31) and at the second inequality by the evident estimate for $\#S_0^{(p)}$ based on (28). Finally, finiteness of the sum over $S_2$ follows from (29) and Lemma 5.1(ii). $\square$

We note in passing that in the proof of Proposition 6.1, we over-estimated the left-hand side of (27) by requiring in (30) that $|\delta\mathbf{j}(e)| \leq b$ only for $e \in E_0$, rather than for all $e \in \Pi_0$.

PROPOSITION 6.3. *We continue under the hypotheses of the preceding proposition, and now make the further assumption that $\#\Pi_1 \vee \Pi = 1$. Then: we have*

$$\sum_{\substack{\mathbf{j}:\,(p,2k)\text{-word s.t. } \mathbf{j} \\ \text{is } \Pi_1\text{-measurable}}} (1 - B(\mathbf{j}))C_\Pi(\mathbf{j}) = o_{p\to\infty}(p). \tag{32}$$

PROOF. We continue in the graph-theoretical setup of the proof of the preceding proposition. But now, under our additional hypothesis that $\#\Pi_1 \vee \Pi = 1$, the set $E_0$ is empty, and hence the set $S_0^{(p)}$ is now simply the set of constant functions on $\{1, \ldots, 2k\}$ taking values in the set $\{1, \ldots, p\}$. Fix $\epsilon > 0$ arbitrarily and then choose a finite set $F \subset S_2$ such that $\sum_{\mathbf{j} \in S_2 \setminus F} |C_\Pi(\mathbf{j})| < \epsilon$. Let

$$N = \max\{|\mathbf{j}(\alpha) - \mathbf{j}(\beta)| \mid \alpha, \beta \in \{1, \ldots, 2k\}, \mathbf{j} \in F\}.$$

Let $\mathbf{j}$ be a $\Pi_1$-measurable $(p, 2k)$-word and write $\mathbf{j} = \mathbf{j}_0 + \mathbf{j}_2$ with $\mathbf{j}_0$ a constant function with values in $\{1, \ldots, p\}$ and $\mathbf{j}_2 \in S_2$. If $\mathbf{j}_2 \in F$ then, provided $p$ is large enough to guarantee that $b > N$, we automatically have $B(\mathbf{j}) = 1$. Thus the sum in question is bounded in absolute value by $\epsilon p$ for $p \gg 0$. Since $\epsilon$ is arbitrary, the proposition is proved. □

The proof of the following proposition is immediate from the definitions and therefore omitted.

PROPOSITION 6.4. *Under exactly the same hypotheses as the preceding proposition we have*

$$\lim_{p\to\infty} \frac{1}{p} \sum_{\substack{\mathbf{j}:\,(p,2k)\text{-word s.t. } \mathbf{j} \\ \text{is } \Pi_1\text{-measurable}}} C_\Pi(\mathbf{j}) = \sum_{\mathbf{j} \in \mathbb{Z}^{\Pi_1}/\mathbb{Z}^{\Pi_1 \vee \Pi}} C_\Pi(\mathbf{j}). \tag{33}$$

Lemma 5.1 guarantees that the sum on the right is well defined.

**7. Proof of the law of large numbers.** This section is devoted to the proof of Theorem 2.3. The main point of the proof is summarized by the following result.

PROPOSITION 7.1. *Let Assumptions 2.1 and 2.2 hold. Let $Y = Y^{(p)}$ be as in (4). Let $R_0^{(k)}$ be as in (8). Then we have*

$$\lim_{p\to\infty} p^{-1} \mathbf{E} \operatorname{trace} Y^k = R_0^{(k)} \tag{34}$$

*for every integer $k > 0$.*

From the case $r=2$ of Proposition 8.1, which is proved in the next section, it follows that

$$\lim_{p\to\infty} \text{Var}\left(\frac{1}{p}\text{trace}\,Y^k\right) = 0 \tag{35}$$

for all integers $k>0$. Arguing just as at the end of the proof of [1], Theorem 3.2, one can then deduce Theorem 2.3 from equations (7), (34) and (35). We omit those details. Thus, to finish the proof of Theorem 2.3, we just have to prove Proposition 7.1. (There will be no circularity of reasoning since the proof of Proposition 8.1 does not use Theorem 2.3.)

PROOF OF PROPOSITION 7.1. Back in the setting of Proposition 5.2 with $r=1$ (in which case, $\#\Pi_0 \vee \Pi_1 = 1$), we have

$$\frac{1}{p}\mathbf{E}\,\text{trace}\,Y^k = \sum_{\Pi \in \text{Part}_2(2k)} p^{-1} n^{-k+\#\Pi_0 \vee \Pi} \sum_{\substack{\mathbf{j}:\,(p,2k)\text{-word s.t. }\mathbf{j}\\ \text{is }\Pi_1\text{-measurable}}} B(\mathbf{j}) C_\Pi(\mathbf{j}).$$

For fixed $\Pi \in \text{Part}_2(2k)$ the contribution to the total sum is

$$O\left(n^{-1-k+\#\Pi_0 \vee \Pi + \#\Pi_1 \vee \Pi}\left(\frac{b}{n}\right)^{-1+\#\Pi_1 \vee \Pi}\right)$$

by Proposition 6.1. Thus, in view of Proposition 3.1, specifically estimate (13), in order to evaluate the limit in question, we can throw away all terms, save those associated to $\Pi = \Pi_0$. We therefore have

$$\lim_{p\to\infty} \frac{1}{p}\mathbf{E}\,\text{trace}\,Y^k = \sum_{\mathbf{j} \in \mathbb{Z}^{\Pi_1}/\mathbb{Z}^{\Pi_0 \vee \Pi_1}} C_{\Pi_0}(\mathbf{j}) \tag{36}$$

by Propositions 6.3 and 6.4. Recalling that $R(j-i) = \mathbf{C}(Z_i, Z_j)$, and writing

$$\mathbf{j} = (j_1, j_2, j_2, \ldots, j_k, j_k, j_1)$$

we have

$$C_{\Pi_0}(\mathbf{j}) = R(j_2 - j_1) \cdots R(j_k - j_{k-1}) R(j_1 - j_k),$$

and hence

$$\sum_{\mathbf{j} \in \mathbb{Z}^{\Pi_1}/\mathbb{Z}^{\Pi_0 \vee \Pi_1}} C_{\Pi_0}(\mathbf{j}) = \sum_{j_2,\ldots,j_k \in \mathbb{Z}} R(j_2 - j_1) \cdots R(j_k - j_{k-1}) R(j_1 - j_k) = R_0^{(k)}$$

for any fixed $j_1 \in \mathbb{Z}$. The proof of (34) is complete. □

**8. Proof of the central limit theorem.** This section is devoted to the proof of Theorem 2.4. The main point of the proof is summarized by the following proposition.

PROPOSITION 8.1. *Let Assumptions* 2.1 *and* 2.2 *hold. Let* $Y = Y^{(p)}$ *be as in* (4). *Let* $Q_{ij}$ *and* $R_i^{(m)}$ *be as in* (8). *Then for each integer* $r \geq 2$, *and all positive integers* $k_1, \ldots, k_r$, *we have*

$$\lim_{p \to \infty} \left(\frac{n}{p}\right)^{r/2} \mathbf{C}(\operatorname{trace} Y^{k_1}, \ldots, \operatorname{trace} Y^{k_r})$$
$$= \begin{cases} 0, & \text{if } r > 2, \\ k_1 k_2 \left( 2R_0^{(k_1+k_2)} + \sum_{i,j} R_i^{(k_1-1)} Q_{ij} R_j^{(k_2-1)} \right), & \text{if } r = 2. \end{cases}$$

In view of Lemma 4.2, in order to finish the proof of Theorem 2.4 by the method of moments, we just have to prove Proposition 8.1.

PROOF OF PROPOSITION 8.1. Back in the setting of Proposition 5.2, this time assuming $r \geq 2$, we have

$$\left(\frac{n}{p}\right)^{r/2} \mathbf{C}(\operatorname{trace} Y^{k_1}, \ldots, \operatorname{trace} Y^{k_r})$$
$$= \sum_{\substack{\Pi \in \operatorname{Part}_2(2k) \\ \text{s.t. } \#\Pi_0 \vee \Pi_1 \vee \Pi = 1}} p^{-r/2} n^{r/2 - k + \#\Pi_0 \vee \Pi} \sum_{\substack{\mathbf{j}: (p, 2k)\text{-word s.t. } \mathbf{j} \\ \text{is } \Pi_1\text{-measurable}}} B(\mathbf{j}) C_\Pi(\mathbf{j}),$$

(37)

and for fixed $\Pi$ the contribution to the total sum is

$$O\left( p^{1-r/2} n^{r/2 - k - 1 + \#\Pi_0 \vee \Pi + \#\Pi_1 \vee \Pi} \left(\frac{b}{n}\right)^{-1 + \#\Pi_1 \vee \Pi} \right)$$

by Proposition 6.1. In view of Proposition 3.1, specifically estimate (14), we are already done in the case $r > 2$.

For the rest of the proof assume $r = 2$. By the estimate immediately above many terms can be dropped from the right-hand side of the sum (37) without changing the limit as $p \to \infty$. The terms remaining can be analyzed by means of Propositions 3.1, 6.3 and 6.4. We thus obtain the formula

(38) $$\lim_{p \to \infty} \frac{n}{p} \mathbf{C}(\operatorname{trace} Y^{k_1}, \operatorname{trace} Y^{k_2}) = \sum_{\substack{\Pi \in \operatorname{Part}_2(2k) \\ \text{s.t. } \#\Pi_1 \vee \Pi = 1 \\ \text{and } \#\Pi_0 \vee \Pi = k - 1}} K(\Pi)$$

where

$$K(\Pi) = \sum_{\mathbf{j} \in \mathbb{Z}^{\Pi_1} / \mathbb{Z}^{\Pi_1 \vee \Pi}} C_\Pi(\mathbf{j}).$$

It remains only to classify the $\Pi$'s appearing on the right-hand side of (38) and for each to evaluate $K(\Pi)$.

We turn to the classification of $\Pi$ appearing on the right-hand side of (38). Recall that in the setup of Proposition 5.2 with $r = 2$, we have

$$\Pi_0 = \{\{1,2\},\ldots,\{2k-1,2k\}\},$$
$$\Pi_1 = \{\{2,3\},\ldots,\{2k_1,1\},\{2k_1+2,2k_1+3\},\ldots,\{2k,2k_1+1\}\}.$$

The conditions

$$\#\Pi_0 \vee \Pi = k-1, \qquad \#\Pi_1 \vee \Pi = 1$$

dictate that we must have

$$(\Pi_0 \vee \Pi) \setminus \Pi_0 = \{A \cup A'\}, \qquad \Pi_0 \setminus (\Pi_0 \vee \Pi) = \{A, A'\}$$

for some $A, A' \in \Pi_0$ with

$$A \subset \{1,\ldots,2k_1\}, \qquad A' \subset \{2k_1+1,\ldots,2k\}.$$

There are exactly $k_1 k_2$ ways of choosing such $A$ and $A'$, and for each such choice, there are exactly three possibilities for $\Pi$, two of which are perfect matchings and one which has all parts of size 2 except for one part of size 4. That is, either

$$(39) \qquad \Pi = (\Pi_0 \setminus \{A, A'\}) \cup \{\{\min A, \min A'\}, \{\max A, \max A'\}\}$$

or

$$(40) \qquad \Pi = (\Pi_0 \setminus \{A, A'\}) \cup \{\{\min A, \max A'\}, \{\max A, \min A'\}\}$$

or

$$(41) \qquad \Pi = (\Pi_0 \setminus \{A, A'\}) \cup \{A \cup A'\}.$$

Thus we have enumerated all possible $\Pi$'s appearing on the right-hand side of formula (38). We remark that Figure 1 depicts examples of $\Pi$ falling into patterns (41), (39), respectively.

We turn to the evaluation of $K(\Pi)$ in the cases (39), (40). In these cases, simply because $\#\Pi \vee \Pi_1 = 1$ and $\Pi$ is a perfect matching, it is possible to choose a permutation $\sigma$ of $\{1,\ldots,2k\}$ such that

$$\Pi_1 = \{\{\sigma(2),\sigma(3)\},\ldots,\{\sigma(2k),\sigma(1)\}\},$$
$$\Pi = \{\{\sigma(1),\sigma(2)\},\ldots,\{\sigma(2k-1),\sigma(2k)\}\},$$

and so we find in these cases that

$$(42) \qquad\qquad\qquad K(\Pi) = R_0^{(k)}$$

by a repetition of the calculation done at the end of the proof of Proposition 7.1.

We turn finally to the evaluation of $K(\Pi)$ in the case (41). In this case there is enough symmetry to guarantee that $K(\Pi)$ does not depend on $A$ and $A'$. We may therefore assume without loss of generality that

$$A = \{2k_1 - 1, 2k_1\}, \qquad A' = \{2k_1 + 1, 2k_1 + 2\}$$

in order to evaluate $K(\Pi)$. To compress notation we write

$$C_{j_1 j_2 j_3 j_4} = \mathbf{C}(Z_{j_1}, Z_{j_2}, Z_{j_3}, Z_{j_4}), \qquad R_{ij}^{(m)} = R_{j-i}^{(m)}, R_{ij} = R_{ij}^{(1)}.$$

Assume temporarily that $k_1, k_2 > 1$. Since $R_{ij} = \mathbf{C}(Z_i, Z_j)$ we then have for any fixed $j_1 \in \mathbb{Z}$ that

$$K(\Pi) = \sum_{j_2,\ldots,j_k \in \mathbb{Z}} R_{j_1 j_2} \cdots R_{j_{k_1-1} j_{k_1}} C_{j_{k_1} j_1 j_{k_1+1} j_{k_1+2}} R_{j_{k_1+2} j_{k_1+3}} \cdots R_{j_k j_{k_1+1}}$$

and hence after summing over "interior" indices we have

(43) $\quad K(\Pi) = \displaystyle\sum_{j_2, j_3, j_4 \in \mathbb{Z}} R_{j_1 j_2}^{(k_1-1)} C_{j_2 j_1 j_3 j_4} R_{j_4 j_3}^{(k_2-1)} = \sum_{i,j} R_i^{(k_1-1)} Q_{ij} R_j^{(k_2-1)}.$

One can then easily check by separate arguments that (43) remains valid when $k_1$ or $k_2$ or both take the value 1.

Together (38)–(43) complete the proof. $\square$

## 9. Concluding comments.

1. We have presented a combinatorial approach to the study of limits for the spectrum of regularized covariance matrices. We have chosen to present the technique in the simplest possible setting, that is, the stationary setup with good a-priori estimates on the moments of the individual entries. Some directions for generalization of this setup are to allow nonstationary sequences with covariances, as in [5], or to allow for perturbations of the stationary setup, as in [3], or to relax the moment conditions of Assumption 2.2. In these more general situations, especially in the context of the LLN, the techniques we have presented here, plus standard approximation techniques, are likely to yield results. But to keep focused we do not study these here.

2. A natural question is whether our approach applies also to the study of centered empirical covariances, that is matrices $\tilde{Y} = \tilde{Y}^{(p)}$ with entries

$$\tilde{Y}(i,j) = B_{ij}(X - \tilde{X})^T (X - \tilde{X})_{ij},$$

where $\tilde{X}_{ij} = n^{-1} \sum_{k=1}^n X_{kj} = n^{-3/2} \sum_{k=1}^n Z_j^{(k)}$. To a limited extent it does, as we now explain. Note that with $(\Delta^{(p)})_{ij} = n B_{ij}^{(p)} \tilde{X}_{1i} \tilde{X}_{1j} = B_{ij}^{(p)} (n^{-1} \sum_{k=1}^n Z_i^{(k)}) \times$ $(n^{-1}\sum_{k=1}^n Z_j^{(k)})$ we have $\tilde{Y} = Y - \Delta^{(p)}$. In contrast to the non-banded version, the perturbation $\Delta = \Delta^{(p)}$ is not of rank one. However, for some deterministic constants $C_1$, $C_2$,

$$\mathbf{E}\sum_{i,j=1}^p \Delta(i,j)^2 \leq C_1 b \sum_{j=1}^p \mathbf{E}\left(\frac{1}{n}\sum_{i=1}^n Z_j^{(i)}\right)^4 \leq \frac{C_2 bp}{n^2}.$$

In particular,

$$\sqrt{\sum_{i,j=1}^p \Delta(i,j)^2} = O_P(\sqrt{bp}/n),$$

where for a sequence of positive random variables $W_p$ and deterministic positive sequence $g_p$ we say that $W_p = O_P(g_p)$ if $W_p/g_p$ is a tight sequence as $p \to \infty$. Letting $\lambda_1 \leq \lambda_2 \leq \cdots \leq \lambda_p$ denote the (ordered) collection of eigenvalues of $Y^{(p)}$ and $\tilde{\lambda}_1 \leq \tilde{\lambda}_2 \leq \cdots \leq \tilde{\lambda}_p$ that of $\tilde{Y}^{(p)}$, we conclude that $\sum_{i=1}^p |\tilde{\lambda}_i - \lambda_i|^2 = O_P(bp/n^2)$. Let $\mathrm{Lip}_1$ denote the collection of deterministic Lipschitz functions on $\mathbb{R}$ bounded by 1 and with Lipschitz constant 1. Since

$$\sqrt{\frac{n}{p}} \sup_{f \in \mathrm{Lip}_1} \sum_{i=1}^p |f(\lambda_i) - f(\tilde{\lambda}_i)| \leq \sqrt{n}\sqrt{\sum_{i=1}^p |\lambda_i - \tilde{\lambda}_i|^2} = O_P(\sqrt{bp/n}),$$

it follows from Theorem 2.3 that the conclusion of that theorem remain true if we substitute $\tilde{L}$, the empirical measure of the eigenvalues of $\tilde{Y}^{(p)}$, for $L$. But this line of reasoning is not powerful enough to yield the conclusion of Theorem 2.4 with $\tilde{Y}$ replacing $Y$, unless one imposes additional conditions that, in particular, imply $bp/n \to_{p\to\infty} 0$. Thus it appears that our basic method itself requires some significant modification to yield a CLT in the regime of Assumption 2.1 in the centered case. The problem is open.

3. We emphasize that unlike the results in [3], we do not deal at all with the distance (in operator norm, or otherwise) between the banded empirical covariance matrix $Y$, and the covariance matrix of the process $\{Z_j\}$.

4. A natural question arising from the central limit theorem (Theorem 2.4) is whether one can obtain an approximation for $E\operatorname{trace} Y^k$ with $o_{p\to\infty}(1)$ error. We recall that in the context of classical Wishart matrices, compact formulas for these quantities can be written down; see [1] and references therein. A similar attempt to provide such formulas here seems to run into many subcases, depending on the relations between the parameters $p, n, b$, and on the convergence rate in the summability condition. We were unsuccessful in finding compactly expressible results. We thus omit this topic entirely.

5. We finally mention a combinatorial question arising from Proposition 3.1. In the setting of that proposition, it can be shown that for perfect matchings $\Pi$ the estimate

$$(44) \qquad \#\Pi_0 \vee \Pi + \#\Pi_1 \vee \Pi \leq k + 2 - r$$

holds and is sharp. But (44) is too strong to hold in general, as is shown by the example

$$\Pi_0 = \{\{1,2\}, \{3,4\}, \{5,6\}, \{7,8\}, \{9,10\}, \{11,12\}\},$$
$$\Pi_1 = \{\{2,3\}, \{1,4\}, \{5,6\}, \{7,8\}, \{9,10\}, \{11,12\}\},$$
$$\Pi = \{\{1,5,6\}, \{2,7,8\}, \{3,9,10\}, \{4,11,12\}\}$$

for which

$$\#\Pi_0 \vee \Pi = \#\Pi_1 \vee \Pi = 2, \qquad k = 6, r = 5,$$

and the same example leaves open the possibility that (14) is too weak. How then can one sharpen (14)? The problem is open.

**Acknowledgments.** We thank Peter Bickel and Noureddine El-Karoui for many discussions that led us to consider the problem discussed in this paper. O. Z. thanks SAMSI for facilitating this exchange.

SCHOOL OF MATHEMATICS
UNIVERSITY OF MINNESOTA
206 CHURCH ST. SE
MINNEAPOLIS, MINNEAPOLIS 55455
USA
E-MAIL: gwanders@math.umn.edu
zeitouni@math.umn.edu